\documentclass[12pt]{amsart}
\usepackage{amscd,amssymb}
\usepackage[arrow,matrix]{xy}
\usepackage[plainpages,backref,urlcolor=blue]{hyperref}

\theoremstyle{plain}

\newtheorem{lem}[subsection]{Lemma}
\newtheorem{prop}[subsection]{Proposition}
\newtheorem{cor}[subsection]{Corollary}

\theoremstyle{definition}
\newtheorem{rk}[subsection]{Remark}

\numberwithin{equation}{section}
\newcommand{\C}{\mathbb{C}}

\newcommand{\PP}{\mathbb{P}}

\begin{document} 

\title [Some criteria to check if a projective hypersurfaces is smooth or singular]
{Some criteria to check if a projective hypersurfaces is singular or smooth} 

\author[Gabriel Sticlaru]{Gabriel Sticlaru}
\address{Faculty of Mathematics and Informatics,
Ovidius University,
Romania}
\email{gabrielsticlaru@yahoo.com }



\begin{abstract}
In this paper we present some properties for projective hypersurfaces, smooth and singular, to be criteria for identification. To make the decision with these criteria, we have included procedures written in Singular language.  
 
\end{abstract}

\maketitle
 
\section{Introduction}  

Let $S=\C[x_0,...,x_n]$ be the graded ring of polynomials in $x_0,,...,x_n$ with complex coefficients and denote by $S_d$ the vector space of homogeneous polynomials in $S$ of degree $d$. 
For any polynomial $f \in S_d$ we define the {\it Jacobian ideal} $J_f \subset S$ as the ideal spanned by the partial derivatives $f_0,...,f_n$ of $f$ with respect to $x_0,...,x_n$.

The Hilbert-Poincar\'{e} series of a graded $S$-module $M$ of finite type is defined by 
\begin{equation} 
\label{eq1}
HP(M)(t)= \sum_{k\geq 0} \dim M_k\cdot t^k 
\end{equation} 
and it is known, to be a rational function of the form 
\begin{equation} 
\label{eq2}
HP(M)(t)=\frac{P(M)(t)}{(1-t)^{n+1}}.
\end{equation}

For any polynomial $f \in S_d$ we define the hypersurface $V(f)$ given by $f=0$ in $\PP^n$ and the corresponding graded {\it Milnor} (or {\it Jacobian}) {\it algebra} by
\begin{equation} 
\label{eq3}
M(f)=S/J_f.
\end{equation}

Smooth hypersurfaces $V(f_s)$ of degree $d$ have all the same  Hilbert-Poincar\'{e} series.

\begin{prop} \label{prop1}

The following statements are equivalent:

\noindent (i) The hypersurface V(f) is smooth.

\noindent (ii) The Hilbert-Poincar\'{e} series is 
$$HP(M(f))(t)=\frac{(1-t^{d-1})^{n+1}}{(1-t)^{n+1}}=(1+t+t^2+\ldots + t^{d-2} )^{n+1}$$
\end{prop}

For a proof see, for instance,   \cite{D1}, p. 109.

As soon as the hypersurface $V(f)$ acquires some singularities, the series $HP(M(f))$ is an infinite sum.

Beyond the smooth case, there is just one general situation in which we know explicit formulas for $HP(M(f))$. 
For polynomials $f$ such that the saturated Jacobian ideal $\widehat J_f$ of the Jacobian ideal $J_f$ is a complete intersection ideal of multidegree $(d_1,...,d_n)$, then it is known that
\begin{equation} 
\label{eq4}
HP(M(f))(t)=\frac{(1-t^{d-1})^{n+1}+t^{(n+1)(d-1)-\sum d_i} (1-t^{d_1}) \cdots (1-t^{d_n})}{(1-t)^{n+1}},
\end{equation}
see Proposition 4, in  \cite{D3}.

For some curves or surfaces, we can find effective singularities and classify.

It is generally difficult to find all the singularities. In this paper we will check if there is at least one singular point and possibly how many they are.
It is a very hard work to obtain manually these informations, even if we consider curves and surfaces with low degree. 
In the last part, we present procedures  in Singular languages.

\section{Smooth hypersurfaces}  \label{smooth}

Smooth hypersurfaces have the same  Hilbert-Poincar\'{e} series associated to Fermat type

$F(t)=(1+t+t^2+\ldots + t^{d-2} )^{n+1} =\sum_{k=0}^{k=T}a_{k}t^{k}$, where $T=(n+1)(d-2)$.

To the best of our knowledge, the only results are for $n=2:$
 \[ a_k = 
     \begin{cases}
	    \binom{k+2}{2}, 0 \leq k \leq d-2 \\
	   \binom{k+2}{2}-3\binom{k+3-d}{2}, d-1 \leq k \leq T/2 \\
      a_{T-k},  \frac{T}{2} < k \leq T.
	   \end{cases}
\] 

In this section we show explicit formulas for the coefficients $a_k$. To compute these coefficients, we need additional informations.

\begin{lem}
The number of solutions in positive integers ($b_{i}\geq 1,k \geq n+1$) of
the equation $b_{0}+b_{1}+\ldots +b_{n}=k$ is ${k-1\choose n}$.
\end{lem}

\proof
Let be the sequences $111\ldots \ 111$ \ with $k$ bits and $k-1$ spaces. If
we keep only $n$ spaces (and all other will be deleted), we have $n+1$
groups of successive bits. There are obviously ${k-1\choose n}$ ways of doing this.
\endproof

\begin{cor}
The number of solutions in nonnegative integers ($b_{i}\geq 0,k\geq 0)$ of
the equation $b_{0}+b_{1}+\ldots +b_{n}=k$ is ${n+k\choose n}$.
\end{cor}

\proof
The transformation $b_{i}=c_{i}-1,i=0,1,\ldots ,n$ yields an equation \ $%
c_{0}+c_{1}+\ldots +c_{n}=n+k+1$ with all $c_{i}\geq 1$. The number of
solutions of this equivalent equation is ${n+k\choose n}$.
\endproof

\begin{rk}
If $k\in 
\mathbb{N}
,$ $k\leq d-2,$ the number of solutions in nonnegative integers, $0\leq
b_{i}\leq d-2$ of the equation $b_{0}+b_{1}+\ldots +b_{n}=k$ is ${n+k\choose n}$. 
We note that condition $b_{i}\leq d-2$ is not necessary.
\end{rk}

\begin{rk}
If $k\in 
\mathbb{N}
,$ $k>(d-2)(n+1),$ then the equation $b_{0}+b_{1}+\ldots +b_{n}=k$ with $%
0\leq b_{i}\leq d-2,$ has no solutions.
\end{rk}

\begin{lem}
The number of solutions in nonnegative integers ($b_{i}\geq 0)$ of the
equation $b_{0}+b_{1}+\ldots +b_{n}=k$ with $m \geq d-1$ variables  and the
other $n+1-m \leq d-2$ variables  is ${n+k-m(d-1)\choose n}$.
\end{lem}

\proof
We can permute the variables, and first take $m\geq d-1.$ The
transformation $b_{i}=c_{i}+(d-1),i=0,1,\ldots ,m-1$ and $b_{i}=c_{i}$ for $%
i=m,\ldots ,n$, yields an equation $c_{0}+c_{1}+\ldots +b_{n}=k-m(d-1).$
Because all $c_{i}\geq 0$ we have ${n+k-m(d-1)\choose n}$ solutions.
\endproof

Now we can prove the main result.

\begin{prop} \label{myprop1}
Let $f=F_{d}(x_{0},x_{1}\ldots x_{n})=x_{0}^{d}+x_{1}^{d}+\cdots +x_{n}^{d}$ 
$\in 
\mathbb{C}
\lbrack x_{0},x_{1}\ldots x_{n}]$. Then $HP(M(f))=\sum_{k\geq 0}a_{k}t^{k}=F(t)$
with $F(t)$ a polynomial of degree $T=(d-2)(n+1)$. In addition

\begin{itemize}
\item for $k=0,\ldots ,d-2$, $a_{k}={n+k\choose n}$.

\item for $k\geq d-1$, $a_{k}={n+k\choose n}+%
\sum_{m=1}^{q}(-1)^{m}{n+1\choose m}\cdot {n+k-m(d-1)\choose n}$ where $q=\left[ \frac{k%
}{d-1}\right] $
\end{itemize}

\end{prop}

\proof
Because $a_{k}=\dim_{%
\mathbb{C}
}(M(f)_{k})$ then $a_{k}=Card(M_{k}$) where $M_{k}$ is the set of solutions
for the equation: $b_{0}+b_{1}+\ldots +b_{n}=k$ with $0\leq b_{i}\leq d-2$
and $0\leq k\leq (n+1)(d-2)$. Obviously $a_{k}={n+k\choose n}$ for $k\in 
\mathbb{N}
,k\leq d-2$. For the second part, $k\geq d-1$, is necessary to involve
principle of inclusion-exclusion. If $B$ is the set of the solutions of the
equation, without the restrictions $b_{i}\leq d-2$, then $B$ has ${n+k\choose n}$
elements. Because $k\geq d-1$, we have some i with $b_{i}\geq d-1.$ Let $%
B_{i}$ be the set of this solutions with $b_{i}\geq d-1$ then $B_{1},\ldots
,B_{n}$ have the same cardinal, ${n+k-(d-1)\choose n}$  and the intersection of $m$
sets from $B_{0},B_{1},\ldots ,B_{n}$ has ${n+k-m(d-1)\choose n}$ elements. If we
have $m$ variables with $b_{i}\geq d-1$, then $m \leq \frac{k}{d-1}$ (othewise, 
$b_{0}+b_{1}+\ldots +b_{n}>k$), so $1\leq m\leq \left[ \frac{k}{d-1}\right] $%
. Hence $M_{k}$ is the difference of $B$ and the union of $B_{0},B_{1},\ldots
,B_{n}$, based on inclusion-exclusion principle, we obtain:
$$a_{k}=Card(M_{k})= {n+k\choose n}+\sum_{m=1}^{q}(-1)^{m}{n+1\choose m}\cdot
{n+k-m(d-1)\choose n}.$$
\endproof

\noindent 
\textbf { Formulas for some cases }

We need to compute the $a_k$ coefficients only for  $ 0\leq k\leq \frac{T}{2}.$ because
$a_{k}=a_{T-k}$  for  $k> \frac{T}{2}.$

\noindent 
$\bullet$ Case $n=3$  \\
$a_{k}={k+3\choose 3}$ for $0\leq k\leq d-2$.\\
$a_{k}={k+3\choose 3}-4 \cdot {k+4-d \choose 3}$ for $d-1\leq k\leq \frac{T}{2}=2d-4< 2d-3$.\\
\medskip 

\noindent 
$\bullet$ Case $n=4$  \\
$a_{k}={k+4\choose 4}$ for $0\leq k\leq d-2$. \\
$a_{k}={k+4\choose 4}- 5 \cdot {k+5-d \choose 4}$ for $d-1\leq k\leq 2d-3$.\\
$a_{k}={k+4\choose 4}- 5 \cdot {k+5-d \choose 4}+ 10 \cdot {k+6-2d \choose 4}$ for $2d-2\leq k\leq \frac{T}{2}<3d-4$.\\

\section{Some criteria and examples}

In this section we present some criteria to check if the projective hypersurfaces is singular or smooth.

\subsection{Critical points}

For some curves or surfaces, we can find effective singularities and classify.
For low degree and simple polynomials, we can find the critical points by solving the system of equations: 
\begin{equation} 
\label{eq_sing}
\frac{\partial f}{\partial x_0} =0, \dots, \frac{\partial f}{\partial x_n} =0. 
\end{equation}

Because $f$ is homogenous of degree $r$, by Euler's formula we have: $x_0 \frac{\partial f}{\partial x_0}+ \dots + 
x_n \frac{\partial f}{\partial x_n} = r f.$ Obviously, if all partial derivatives of polynomial $f$ vanish at $p$, then the polynomial $f$ vanishes at $p$ too. Therefore, for finding singular points it is enough to get critical points and classify.

For example, in $\PP^2$, these points $(0:0:1), (1:i:0), (1:-i:0)$ are nodes (type $A_1$) for Lemniscate of Bernoulli: $(x^2+y^2)^2-2(x^2-y^2)z^2=0$   and cusps (type $A_2$) for Cardioid: $(x^2+y^2+xz)^2 - (x^2+y^2)z^2=0. $ 

\textbf{Remark}

It is sufficient to find at least one critical point.
For projective surfaces in $\PP^3$, $f=(x^4+w^4)^p+(y^p+z^p)^4=0$ it is easy to check that $(1:0:0:i)$ is critical
point, for arbitrary $p$. 

\subsection{Hilbert-Poincar\'{e} series}

The Hilbert-Poincar\'{e} series is computed with two methods: combinatorial and based on a free resolution.
For other examples see \cite{DSt2}, \cite{DSt3}.

To compute the Hilbert-Poincar\'{e} series, first recall that  the quotient rings $S/I$ and $S/LT(I)$ have the same series, for any monomial ordering, where LT(I) is the ideal of leading terms of the ideal I. 

The projective equation of the Cissoid of Diocles is given by $x^3+xy^2-y^2z=0$.

The leading ideal of the jacobian $J_f $ is: $ LI=(y^2, xy, x^2)$.

For the graded Milnor algebra, $M=\oplus_{k\geq0}M_{k} $ we show the bases for the  homogeneous components: $ M_0= \left\{  1 \right\}, M_1= \left\{   z,y,x  \right\},M_2= \left\{   z^2,yz,xz  \right\},
M_3= \left(  z^3,xz^2  \right\},\\
M_4= \left(  z^4,xz^3  \right\}$ and  $M_k= \left\{   z^k,xz^{k-1}  \right\}$ for all $k\geq 3$.

Finally, if we count the number of monomials in each homogeneous bases, we find the Hilbert-Poincar\'{e} series 
$ S(t)=1+3t+3t^2+2(t^3+\ldots $ 

The projective equation of the curve called Kappa is $ f=(x^2+y^2)y^2-x^2z^2=0$.

Here is the minimal graded free resolution of Milnor algebra M:
\begin{equation} 
\label{res}
0 \to R_3 {\rightarrow} R_2 {\rightarrow} R_1 {\rightarrow} R_0 \to M \to 0
\end{equation}
where $R_0=S$, $R_1=S^3(-3)$, $R_2=S(-5) \oplus S^3(-6)$ and $R_3= S^2(-7)$.

To get the formulas for the Hilbert-Poincar\'{e} series, we start with the resolution \eqref{res} and get
$HP(M)(t)=HP(R_0)(t)-HP(R_1)(t)+HP(R_2)(t)-HP(R_3)(t).$

Then we use the well-known formulas $HP(N_1\oplus N_2)(t)=HP(N_1)(t)+HP( N_2)(t)$, $HP(N^p(-q))(t)=pt^qHP(N)(t)$,
$HP(S)(t)=\frac{1}{(1-t)^3}$ and we obtain:

$HP(M)(t)=\frac{1-3t^3+t^5+3t^6-2t^7}{(1-t)^3}=\frac{1+2t+3t^2+t^3-t^4-2t^5}{1-t}=1+3t+6t^{2}+7t^{3}+6t^{4}+4(t^{5}+\ldots $

\textbf{Remark}

It is sufficient to find just one coefficient which is different in the two series to infer that $V(f)$ is singular.

For smooth curve of degree three, the Hilbert-Poincar\'{e} series is polynomial, $F_3(t)=(1+t)^3=1+3t+3t^2+t^3.$
For the Cissoid of Diocles curve, $M_3= \left(  z^3,xz^2  \right\}$ with $dim(M_3)=2$ so the Cissoid is singular.

Another criteria is based on degree $T=(n+1)(d-2)$ of the polynomial for the smooth case. If  $dim(M_{T+1})>0$ the hypersurfaces is singular. 

\subsection{Complete intersection}

We present projective hypersurfaces $V(f)$ for which the saturated Jacobian ideal $\widehat J_f$ is a complete intersection of multidegree $(d_1,...,d_n)$, so these $V(f)$ are singular and their Hilbert-Poincar\'{e} series $HP(M(f))(t)$ are infinite. 
We display the type of complete intersection $(d_1,...,d_n),$ the generators for the saturated Jacobian ideal $\widehat J_f$ and the Hilbert-Poincar\'{e} series.

\medskip 

\noindent 
$\bullet$ singular curves in $\PP^2$  \\
\medskip 
  $ V(f): f=x^4+y^4+xyz^2=0$ with $\widehat J_f$ of type (1,1), and two generators: $(x,y)$. 
  
$HP(M(f))(t)=\frac{(1-t^3)^3+t^7 (1-t)(1-t)}{(1-t)^3}
=1+3t+6t^{2}+7t^{3}+6t^{4}+3t^{5}+(t^{6}+t^{7}+\dots $

(this curve has a single node).

\medskip 
\noindent 
$ V(f): f= x^2y^2+xz^3+yz^3=0$ with $\widehat J_f$ of type (2,2), and two generators: $(xy,z^2)$.

$HP(M(f))(t)=\frac{(1-t^3)^3+t^5 (1-t^2)(1-t^2)}{(1-t)^3}=
1+3t+6t^{2}+7t^{3}+6t^{4}+4(t^{5}+t^{6}+\dots $

(this curve has 2 cusps).

\medskip 

\noindent 
$\bullet$ singular surfaces in $\PP^3$  \\
\medskip 
$ V(f): f=(x^2 + y^2 + z^2 + xy+ xz+ yz)w + 2xyz=0$ with $\widehat J_f$ of type (1,1,1), and three generators: 
$(x,y,z)$.
 
$HP(M(f))(t)=\frac{(1-t^2)^4+t^5 (1-t)(1-t)(1-t)}{(1-t)^4}=1+4t+6t^{2}+4t^{3}+(t^{4}+t^{5}+\ldots $

(this surface has a single node).

 \medskip 

\noindent 
$ V(f): f=xzw+(z+w)y^2+x^3+x^2y+xy^2+y^3=0$ with $\widehat J_f$ of type (1,1,2), and three generators: $(x,y,zw)$.

$HP(M(f))(t)=\frac{(1-t^2)^4+t^4 (1-t)(1-t)(1-t^2)}{(1-t)^4}=1+4t+6t^{2}+4t^{3}+2(t^{4}+t^{5}\dots $

(this surface has 2 nodes).

\subsection{Number of singularities}

In this subsection we assume that $V(f)$ has at most isolated singularities.

We consider the following ideals:  $J=J_f$ the jacobian ideal and $I$ the radical of $J$ with $G_1=std(J)$ and 
$G_2=std(I)$ their respective standard Groebner basis, $mult_1$ the degree of $G_1$ and $mult_2$ the degree of $G_2$.

\begin{prop} \label{myprop2}
We have the following statements:
\begin{itemize}
\item If $mult_1 < (d-1)^{n+1}$ the hypersurfaces is singular and $mult_2$ is the number of singularities. 
\item If $mult_1 =(d-1)^{n+1}$ the hypersurfaces is smooth.
\end{itemize}
\end{prop}
\proof

Indeed, when $V(f)$ has isolated singularities, then one has $mult_1=\tau(V(f))$, the total Tjurina number of $V(f)$, which is the sum of the Tjurina numbers $\tau(V(f),p)$ over all singular points $p$ of $V(f)$. This follows from the definition of the multiplicity of a homogeneous ideal and  the equality
$\dim M(f)_k=\tau(V(f))$, for $k$ large, see \cite{CD}. Now, at each singular point we have the inequality
$ \tau(V(f),p) \leq \mu(V(f),p),$
and by taking the sum we get the inequality $ \tau(V(f)) \leq \mu(V(f))$ between the total Tjurina and the total Milnor numbers. Moreover, it is known that $\mu(V(f))\leq (d-1)^n$, see for instance Proposition (3.25), page 90 as well as the discussion on page 161 in \cite{D2}. It follows that in fact, for singular hypersurfaces we have the stronger inequality
$$mult_1 \leq (d-1)^{n}.$$
For the smooth case, the Hilbert-Poincar\'{e} series is polynomial, $F(t)=(1+t+t^2+\ldots + t^{d-2} )^{n+1}$ and we have 
$mult_1= \sum_{k\geq 0} \dim M_k\cdot t^k=F(1)=(d-1)^{n+1}.$

\begin{cor}
\label{corA}
If  the hypersurfaces is singular and $mult_1=mult_2$ them all the singularities are nodes.
\end{cor}

If there are two type of singularities on $V(f)$, say $A_p$ and $A_q$ it is possible to  find the number  $x,y$ of each type of singularities,  by solving the system:
\noindent 
\[ \begin{cases}
	    px+qy=mult_1 \\
	    x+y=mult_2 
	 \end{cases}
\] 

For some case, there are many configuration, for example, if $mult_1 = 5$ and $mult_2 = 3$ we have two configurations 
$A_1+A_2$ and $2A_2+A_3$.
 
In the following table we present some examples for projective surfaces $f=0$ of degree three.

\noindent
\begin{center}
 \begin{tabular} {|l|c|c|l|}
  \hline                
  polynomial $f $   &            $mult_1$ & $mult_2$   & Type   \\  
  \hline
  $ wxz+y^3   				   $  &  $  6  $  & $ 3  $ & $ 3A_2      $ \\
  $ w(xy+xz+yz)+xyz 	   $  &  $  4  $  & $ 4  $ & $ 4A_1      $ \\
  $ wxz+(x+y)y^2   		   $  &  $  5  $  & $ 3  $ & $ A_1+2A_2  $ \\
  $ wxz+(x+z)y^2         $  &  $  5  $  & $ 3  $ & $ 2A_1+A_3  $ \\
  $ wxz+y^2(x+y+z)       $  &  $  4  $  & $ 3  $ & $ 2A_1+A_2  $ \\
  $ wxz+(x+z)(y^2-x^2)   $  &  $  4  $  & $ 2  $ & $ A_1+A_3   $ \\
  $ wxz+y^2z+yx^2        $  &  $  5  $  & $ 2  $ & $ A_1+A_4   $ \\
  $ x^3+y^3+z^3+w^3      $  &  $ 16  $  & $ -  $ & $ smooth    $ \\
  \hline
  
\end{tabular} 
\end{center}

\subsection{Computation of genus}

The genus of a smooth irreducible curve defined by a polynomial of degree
$d$ is given by the formula: $ g_s =\frac{(d-1)(d-2)}{2}$.

The genus of a singular irreducible projective curve $C:f=0$ (i.e. the genus of its smooth model) is given by
$ g(C) = g_s - \delta $, where $\delta$ is the sum of all local delta-invariants $\delta_p$ of the singularities $p \in C$.
If $C$ is smooth, then $\delta$ is 0.

The Hermitian curve  $x^d + y^d = z^d$ and Fermat curve $x^d + y^d + z^d=0$ are smooth in $\PP^2$, 
of genus $\frac{(d-1)(d-2)}{2}$. 

\textbf{Remark}
Every irreducible projective curves with genus=$\frac{(d-1)(d-2)}{2}-1$ is either nodal of type $A_1$ or has exactly one cusp $A_2$.

Several computer algebra packages are able to compute the genus of a plane
curve

In the following table we present some computations with procedure genus() from Singular library normal.lib. 

\noindent
\begin{center}
 \begin{tabular} {|l|c|l|}
  \hline                
  polynomial $f $   &        Genus &  Comments   \\  
  \hline
  $ x^3+xy^2-y^2z`                $  &  $  0  $   & singular ($A_2$) \\
  $ x^3y + y^3z + xz^3            $  &  $  3  $   & smooth (Klein quartic)  \\
  $ x^3yz+y^5+z^5      					  $  &  $  5  $   & singular ($A_1$)  \\
  $ (x^2+y^2)^2-2(x^2-y^2)z^2     $  &  $  0  $   & singular (Lemniscate of Bernoulli)  \\
  $ (x^2+y^2+xz)^2 - (x^2+y^2)z^2 $  &  $  0  $   & singular (Cardioid)  \\
  $ (x^2+y^2)y^2-x^2z^2=0         $  &  $  0  $   & singular (Kappa)  \\
  \hline
\end{tabular} 
\end{center} 
The Lemniscate of Bernoulli, the Cardioid and the Kappa curve are irreducible with genus $=0$ hence these are singular rational (i.e. parameterizable) curves.

\section{ Programs in Singular language} \label{programs}
For mathematical computations, we can use any CAS (Computer Algebra System) software like Mathematica, Matlab or Maple, but for Algebraic Geometry, the best are perhaps Singular, Macaulay2 or CoCoA.

Singular is a package developed at the University of Kaiserslautern, see \cite{DGP} and \cite{GP}. 

\begin{verbatim} 

proc criteria_mult(){  
// compute the number of singularities for projective hypersurfaces
int n=2; // n+1 variables, you can change
ring R=0,(x(0..n)),dp;
poly f;
// you can change polynomial f
f=(x(0)^2+x(1)^2)^2-2*(x(0)^2-x(1)^2)*x(2)^2;  // Lemniscata 
int d=deg(f);
ideal J=jacob(f);
ideal G1=std(J);
ideal G2=std(radical(J));
int mult1=mult(G1);
int mult2=mult(G2);
if(mult1==(d-1)^(n+1)) {print("this curve is smooth!");}
else {print("this curve is singular! ");}
if(mult1==mult2) {print("and nodal! ");}
}

proc criteria_ci(){  
// complete intersection criteria
int n=2; // n+1 variables, you can change
LIB "elim.lib";
ring R=0,(x(0..n)),dp;
poly f;
// you can change polynomial f
f=x(0)^3*x(1)^5+x(1)^8+x(2)^8 ; //(4,7)
int d=deg(f);
ideal J=jacob(f);
ideal I=std(sat(J,maxideal(1))[1]);
if(size(I)!=n){
print("Not Complete intersection:");
return();
}
print("Complete intersection with " + string(n)+" generators:");
print(I);
int i, a;
intvec v;
a=(n+1)*(d-1);
for(i=1;i<=n;i++) {v[i]=deg(I[i]); a=a-v[i];}
ring R=0,t,ds;
poly P, Q;
Q=1;
for(i=1;i<=n;i++) {Q=Q*(1-t^v[i]);}
P=(1-t^(d-1))^(n+1)+ t^a*Q;
print("Hilbert Poincare series P(t)/(1-t)^(n+1) where P(t) is:");
print(P);
}

proc criteria_hp(){  
// compute Hilbert-Poincare series
int n=2; // n+1 variables, you can change
ring R=0,(x(0..n)),dp;
poly f;
// you can change polynomial f
f=(x(0)^2+x(1)^2+x(0)*x(2))^2 - (x(0)^2+x(1)^2)*x(2)^2; //cardioid
int d=deg(f);
ideal J=jacob(f);
ideal G=std(J);
int i, T;
intvec v;
T=(n+1)*(d-2);
for (i=1; i<=T+2;i++) {
v[i]=size(kbase(G,i-1));}
if (v[T+2]==0) {print("Smooth curve!"); return();}
if (v[T+2]>0) {print("The hypersurfaces is singular!");}
if (v[T+2]<>v[T+3]) {print("Not finite isolate singularities!");return();}

ring R=0,t,ds;
poly F, S; 
F=((1-t^(d-1))/(1-t))^(n+1);
S=0;
for(i=1;i<=T+2;i++) {S=S+v[i]*t^(i-1); }
print("Hilbert Poincare series is S(t)= "); S;
print("Hilbert Poincare series smooth (Fermat) case F(t)= "); F;
}

proc criteria_genus (){
// genus, only for projective curve  
LIB "normal.lib";
ring r=0,(x,y,z),dp;
poly f;
// you can change polynomial here
f=(x^2+y^2)^2-2*(x^2-y^2)*z^2; // lemniscata -> genus=0
int d=deg(f);
ideal I=f;
int g=genus(I);
if(g==(d-1)*(d-2)/2) {print("this  curve is smooth !");}
else {print("this curve is singular !");}
} 
\end{verbatim}

\end{document}